\documentclass{article}
\usepackage{amsmath, amsthm}
\usepackage{amssymb}
\usepackage{textcomp}
\newtheoremstyle{theorem}%name
  {10pt}          % space above
  {10pt}  % space below
  {\sl}  % bofy font
  {\parindent}     % ident - empty=no indent,  \parindent= paragraph indent
  {\bf}  % thm head font
  {. }    % punctuation after thm head
  { }    % space after thm head: `` ``=normal \newline=linebreak
  {}     % thm head specification
\theoremstyle{theorem}
\newtheorem{theorem}{Theorem}[section]

\newtheorem{lemma}[theorem]{Lemma}

\newtheoremstyle{defi}%name
  {10pt}          % space above
  {10pt}  % space below
  {\rm}  % bofy font
  {\parindent}     % ident - empty=no indent,  \parindent= paragraph indent
  {\bf}  % thm head font
  {. }    % punctuation after thm head
  { }    % space after thm head: `` ``=normal \newline=linebreak
  {}     % thm head specification
\theoremstyle{defi}

\begin{document}
\title{Higher Order Fractional Variational Optimal Control Problems with Delayed Arguments } \maketitle

\begin{center}
\author\textbf{Fahd Jarad }\\
Department of Mathematics and Computer Sciences, Faculty of Arts
and Sciences, $\c{C}ankaya$ University- 06530, Ankara, Turkey\\

\author\textbf{Thabet Abdeljawad (Maraaba)}\\
Department of Mathematics and Computer Sciences, Faculty of Arts
and Sciences, $\c{C}ankaya$ University- 06530, Ankara, Turkey\\

\author\textbf{Dumitru Baleanu}\footnote{On leave of absence from
Institute of Space Sciences, P.O.BOX, MG-23, R 76900,
Magurele-Bucharest, Romania,E-mails: dumitru@cankaya.edu.tr,
baleanu@venus.nipne.ro}\\
Department of Mathematics and Computer Sciences, Faculty of Arts
and Sciences, $\c{C}ankaya$ University- 06530, Ankara, Turkey\\

\end{center}

\begin{abstract}
This article deals with higher order Caputo fractional variational problems with the presence of delay in the state variables and their integer higher order derivatives.
\end{abstract}

\section{Introduction}
 Lagrangian theories involving higher-order derivatives  appear of great
interest as an imbedding for field theories with fields of subcanonical
dimension,e.g. the Heisenberg nonlinear spinor theory, for which local
interactions are less singular than in the canonical case.

 The calculus of variation has a long history of communications with other fields of mathematics such as geometry and differential equations, and with physics. Recently the calculus of variations has found applications in economics and some branches of engineering such as electrical engineering. Optimal control which is a rapidly expanded field can be regarded as a part of the calculus of variations.

 Recently, the fractional calculus which is as old as the classical calculus has become a candidate in solving problems of complex systems which appear  in various branches of science
\cite{trujillo,samko,podlubny,richard,zaslavsky,machado,mainardi,enrico}.

Several authors found interesting results when they used the fractional calculus in control theory \cite{agrawal4,chen}.

 Experimentally, the use of delay together with  the fractional calculus
may give better results.

  Optimal control problems with time delay in calculus of
variations were discussed in  \cite{rosen}. Variational optimal control problems
within  fractional derivatives were considered in \cite{agragen}.
Fractional variational problems in the presence of delay were
studied in \cite{thab, fahd}.

The aim of this paper is to deal with optimal
control fractional variational problems in the presence of delay in the state variable and its higher order derivatives

The structure of the paper is as follows:

In section 2    the necessary definitions in the  fractional
calculus  used in this manuscript are reviewed. In section 3  the
 unconstrained fractional Euler-Lagrange equations with delay are discussed .  The fractional
control problem is presented in Section 5.

\section{Basic Definitions}

We present in this section some basic definitions related to
fractional derivatives.

\noindent {\em The Left Riemann-Liouville Fractional Integral} and
{\em The Right Riemann-Liouville Fractional Integral} are  defined
respectively by
\begin{equation}
_aI^{\alpha}f(t)=\frac{1}{\Gamma
(\alpha)}\int_a^t(t-\tau)^{\alpha-1}f(\tau)d\tau,
\end{equation}
\begin{equation}
I_b^{\alpha}f(t)=\frac{1}{\Gamma
(\alpha)}\int_t^b(\tau-t)^{\alpha-1}f(\tau)d\tau,
\end{equation}
\noindent where $\alpha >0,~ n-1<\alpha<n$. Here and in the
following $\Gamma(\alpha)$ represents the Gamma function.

 \noindent {\em The Left
Riemann-Liouville Fractional Derivative} is defined by
\begin{equation}
_a{D}^\alpha f(t) = \frac{1}{\Gamma (n-\alpha)} \left(
\frac{d}{dt} \right)^n \int_a^t (t-\tau)^{n-\alpha-1} f(\tau)
d\tau,
\end{equation} %(1)

\noindent {\em The Right Riemann-Liouville Fractional Derivative}
is defined by
\begin{equation}
{D}_b^\alpha f(t) = \frac{1}{\Gamma (n-\alpha)}
\left(-\frac{d}{dt} \right)^n \int_t^b (\tau-t)^{n-\alpha-1}
f(\tau) d\tau,
\end{equation} %(2)

 The fractional derivative of a constant takes the form
\begin{equation}\label{constant}
{}_a{D}^{\alpha}C=\ C\ \frac{(t-a)^{-\alpha}}{\Gamma(1-\alpha)}.
\end{equation}
\noindent and the fractional derivative of a power of t has the
following form
\begin{equation}
{}_a{D}^{\alpha}(t-a)^\beta=\
\frac{\Gamma(\alpha+1)(t-a)^{\beta-\alpha}}{\Gamma(\beta-\alpha+1)},
\end{equation}
for $\beta> -1, \alpha\geq 0$.

The  Caputo's fractional derivatives are defined as
follows:\\

 \noindent {\em The Left Caputo Fractional Derivative}
\begin{equation}
{_a^{C}D^\alpha f(t)} = \frac{1}{\Gamma (n-\alpha)} \int_a^t
(t-\tau)^{n-\alpha-1} \left( \frac{d}{d\tau} \right)^n f(\tau)
d\tau ,
\end{equation} %(3)

\noindent and\\

{\em The Right Caputo Fractional Derivative}
\begin{equation}
{^CD_b^\alpha f(t)} = \frac{1}{\Gamma (n-\alpha)} \int_t^b
(\tau-t)^{n-\alpha-1} \left(-\frac{d}{d\tau} \right)^n f(\tau)
d\tau ,\end{equation} %(4)

\noindent where $\alpha$ represents the order of the derivative
such that $n-1 < \alpha < n$. By definition the Caputo fractional
derivative of a constant is zero.

The Riemann-Liouville fractional derivatives and Caputo fractional
derivatives are connected with each other by the following
relations:
\begin{equation}\label{conn1}
{_a^CD^\alpha f(t)}=~_a{D}^\alpha
f(t)-\sum_{k=0}^{n-1}\frac{f^{(k)}(a)}{\Gamma(k-\alpha+1)}(t-a)^{k-\alpha},
\end{equation}
\begin{equation}\label{conn2}
{^CD_b^\alpha f(t)}={D}_b^\alpha
f(t)-\sum_{k=0}^{n-1}\frac{(-1)^kf^{(k)}(b)}{\Gamma(k-\alpha+1)}(b-t)^{k-\alpha}.
\end{equation}
In \cite{trujillo}, a formula for the fractional
integration by parts on the whole interval $[a,b]$ was given by
the following lemma
\begin{lemma}\label{l1}
 Let $\alpha > 0$, $p,q\geq 1$, and $\frac{1}{p}+
\frac{1}{q} \leq 1+\alpha$ ( $p \ne 1 $ and $q \ne 1 $ in the case
when $\frac{1}{p}+ \frac{1}{q} = 1+\alpha$ )

 (a) If $\varphi \in L_p(a,b)$ and  $\psi \in L_q(a,b)$,
 then
 \begin{equation} \label{e1}
\int_{a}^{b} \varphi(t)(_{a}I^\alpha \psi)(t)dt =\int_{a}^{b}
\psi(t)(I_{b}^\alpha \varphi)(t)dt
 \end{equation}

(b) If $g \in I_{b}^\alpha(L_p)$ and $f \in ~_{a}I^\alpha(L_q) $,
then
\begin{equation} \label{d1}
\int_{a}^{b}g(t) ( {}_{a}{D}^{\alpha} f)(t)dt = \int_{a}^{b}f(t) (
{D}_{b}^\alpha g)(t)dt,
\end{equation}
\noindent where $~_{a}I^\alpha(L_p):=\{f: f=~_{a}I^{\alpha}g, g\in L_p(a,b)\}$ and \\$~I_b^\alpha(L_p):=\{f: f=I_b^{\alpha}g, g\in L_p(a,b)\}$.
\end{lemma}

In  \cite{thab} and \cite{fahd} ,  other formulas for the fractional integration by parts on the subintervals $[a,r]$ and $[r,b]$ were given by the following lemmas

\begin{lemma}\label{l2}
Let $\alpha > 0$, $p,q\geq 1$,~$r \in (a,b)$ and
$\frac{1}{p}+ \frac{1}{q} \leq 1+\alpha$ ( $p \ne 1 $ and $q \ne 1
$ in the case when $\frac{1}{p}+ \frac{1}{q} = 1 + \alpha$ ).

(a)If $\varphi \in L_p(a,b)$ and  $\psi \in L_q(a,b)$,
 then

\begin{equation} \label{f1}
\int_{a}^{r} \varphi(t)(_aI^\alpha \psi)(t)dt =\int_{a}^{r}
\psi(t)(I_r^\alpha \varphi)(t)dt
 \end{equation}

and thus if $g \in ~I_r^\alpha(L_p)$ and $f \in ~_aI^\alpha(L_q)
$, then
\begin{equation} \label{f11}
\int_{a}^r g(t) ( ~_a{D}^{\alpha} f)(t)dt = \int_{a}^rf(t)
({D}_r^\alpha g)(t)dt
\end{equation}

(b)If $\varphi \in L_p(a,b)$ and  $\psi \in L_q(a,b)$,
 then

$$\hspace{-3 cm}\int_r^{b} \varphi(t)(_{a}I^\alpha \psi)(t)dt=
 \int_{r}^{b}
\psi(t)(~I_b^\alpha \varphi)(t)dt$$
\begin{equation}+ \frac{1}
{\Gamma(\alpha)}\int_{a}^r \psi(t) ( \int_r^b \varphi(s)(s-
t)^{\alpha-1}ds )dt
 \end{equation}
and hence if $g \in I_b^\alpha(L_p)$ and $f \in
~_aI^\alpha(L_q)$, then
$$ \hspace{-2.6 cm} \int_r^b g(t) ( ~_aD^{\alpha} f)(t)dt =
\int_r^{b}f(t) ( {D}_b^\alpha g)(t)dt$$
\begin{equation} \label{f11}
-\frac{1} {\Gamma(\alpha)}\int_{a}^r ( ~_a{D}^{\alpha} f)(t) (
\int_r^{b} ( {D}_b^\alpha g)(s)(s- t)^{\alpha-1}ds )dt.
\end{equation}
\noindent That is
$$\hspace{-2.4 cm} \int_r^{b} g(t) (~_a{D}^{\alpha} f)(t)dt =
\int_r^{b}f(t) ({D}_b^\alpha g)(t)dt$$
\begin{equation} \label{f12}
-\frac{1} {\Gamma(\alpha)}\int_{a}^r  f(t) {D}_r^{\alpha}(
\int_r^{b} ( {D}_b^\alpha g)(s)(s- t)^{\alpha-1}ds)dt.
\end{equation}
\end{lemma}

\begin{lemma}\label{l3}
Let $\alpha > 0$, $p,q\geq 1$,~$r \in (a,b)$ and
$\frac{1}{p}+ \frac{1}{q} \leq 1+\alpha$ ( $p \ne 1 $ and $q \ne 1
$ in the case when $\frac{1}{p}+ \frac{1}{q} = 1 + \alpha$ ).

(a)If $\varphi \in L_p(a,b)$ and  $\psi \in L_q(a,b)$,
 then

\begin{equation} \label{f1}
\int_{r}^{b} \varphi(t)(I_b^\alpha \psi)(t)dt =\int_{r}^{b}
\psi(t)(_rI^\alpha \varphi)(t)dt
 \end{equation}

and thus if $g \in ~_rI^\alpha(L_p)$ and $f \in ~I_b^\alpha(L_q)
$, then
\begin{equation} \label{f11}
\int_{r}^b g(t) ( {D}_b^{\alpha} f)(t)dt = \int_{r}^bf(t)
(_r{D}^\alpha g)(t)dt
\end{equation}

(b)If $\varphi \in L_p(a,b)$ and  $\psi \in L_q(a,b)$,
 then

$$\hspace{-3 cm}\int_a^{r} \varphi(t)(_{b}I^\alpha \psi)(t)dt=
 \int_{a}^{r}
\psi(t)(I_{a}^\alpha \varphi)(t)dt$$
\begin{equation}+ \frac{1}
{\Gamma(\alpha)}\int_{r}^b \psi(t) ( \int_a^r \varphi(s)(t-
s)^{\alpha-1}ds )dt
 \end{equation}
and hence if $g \in ~_aI^\alpha(L_p)$ and $f \in
~I_b^\alpha(L_q)$, then
$$ \hspace{-2.6 cm} \int_a^r g(t) ( {D}_b^{\alpha} f)(t)dt =
\int_a^{r}f(t) ( _a{D}^\alpha g)(t)dt$$
\begin{equation} \label{f11}
-\frac{1} {\Gamma(\alpha)}\int_{r}^b ( {D}_b^{\alpha} f)(t) (
\int_a^{r} ( _a{D}^\alpha g)(s)(t- s)^{\alpha-1}ds )dt.
\end{equation}
\noindent That is
$$\hspace{-2.4 cm} \int_a^{r} g(t) ( {D}_b^{\alpha} f)(t)dt =
\int_a^{r}f(t) (  _a{D}^\alpha g)(t)dt$$
\begin{equation} \label{f12}
-\frac{1} {\Gamma(\alpha)}\int_{r}^b  f(t) {}_{r}{D}^{\alpha}(
\int_a^{r} ( _a{D}^\alpha g)(s)(t- s)^{\alpha-1}ds )dt.
\end{equation}
\end{lemma}

\section{ The Unconstrained Caputo Fractional Variation with delay }

Before we consider the fractional control problem, let us consider
the following fractional variational problem with delay

Minimize

 $$J(y)=\int_{a}^{b}L[t,~_a^CD^{\alpha_1}y(t),~_a^CD^{\alpha_2}y(t),...,~_a^CD^{\alpha_n}y(t),$$
$$~^CD_b^{\beta_1}y(t),~^CD_b^{\beta_2}y(t),...
,~^CD_b^{\beta_m}y(t),y(t),y'(t),...,y^{(k)}(t),$$
\begin{eqnarray}\label{23}
y(t-\tau),y'(t-\tau),...,y^{(k)}(t-\tau)]dt,
\end{eqnarray}
such that

$$k-1\le\alpha_{\mbox{max}}<k, ~~\alpha_{\mbox{max}}=\mbox{max}\{\alpha_i,\beta_j\}_{1\le i\le n,1\le j\le m},$$
$$y^{(l)}(b)=c_l,~l=0,1,2,...,k-1,~~y(t)=\phi(t)~ t\in [a-\tau,a],$$
\begin{equation}\label{24}
 \tau>0,~ \tau< b-a,~~\alpha_i,\beta_j\in \mathbb{R}~~\forall ~i=1,2,...,n,~~\forall j=1,2,...,m,
\end{equation}
\noindent where $c_l$ are constant,~$\phi(t)$ is a smooth function and  $L$ is a function with
continuous first and second partial derivatives with respect to
all of its arguments.

If  the above variational problem (\ref{23}) has a minimum at
$y_0(t)$ and $\eta(t)\in \mathbb{R}$ is an admissible function such that
 $\eta(t)\equiv 0$ in the interval $[a-\tau,a]$  then the
function

\begin{eqnarray} \label{25}
\xi (t)=J(y_0+ \epsilon \eta)
\end{eqnarray}
\noindent where $\epsilon \in \mathbb{R}$ admits a minimum at $\epsilon=0$.
Hence

$$\int_a^b[\sum_{i=1}^n \partial_{i+1}L(t)~a^CD^{\alpha_i}\eta(t)+\sum_{j=1}^m \partial_{n+j+1}L(t)^CD_b^{\beta_j}\eta(t)+$$
\begin{equation}\label{26}
\sum_{p=0}^k  \partial_{m+n+p+2}L(t)\eta^{(p)}(t)+\sum_{p=0}^k \partial_{m+n+k+3}L(t)\eta^{(p)}(t-\tau)]dt=0.
\end{equation}
\noindent where $ \partial_i L$ is the partial derivative of $L$ with respect to its $i^{\mbox{th}}$ argument.

On using the connection formulas (\ref{conn1}) and (\ref{conn2}), (\ref{26}) reads

 $$\int_a^b[\sum_{i=1}^n \partial_{i+1}L(t)~aD^{\alpha_i}\eta(t)+\sum_{j=1}^m \partial_{n+j+1}L(t)D_b^{\beta_j}\eta(t)+$$
\begin{equation}\label{27}
\sum_{p=0}^k  \partial_{m+n+p+2}L(t)\eta^{(p)}(t)+\sum_{p=0}^k \partial_{m+n+k+3}L(t)\eta^{(p)}(t-\tau)]dt=0.
 \end{equation}

 Now if one splits the integral, makes the change of variables for $t-\tau$ and uses the fact that $\eta\equiv 0$ in $[a-\tau,a]$, (\ref{27}) becomes
 $$\int_a^{b-\tau}[\sum_{i=1}^n \partial_{i+1}L(t)~aD^{\alpha_i}\eta(t)+\sum_{j=1}^m \partial_{n+j+1}L(t)D_b^{\beta_j}\eta(t)+$$
 $$\sum_{p=0}^k  \partial_{m+n+p+2}L(t)\eta^{(p)}(t)+\sum_{p=0}^k \partial_{m+n+k+3}L(t+\tau)\eta^{(p)}(t)]dt+$$

 $$\int_{b-\tau}^{b}[\sum_{i=1}^n \partial_{i+1}L(t)~aD^{\alpha_i}\eta(t)+\sum_{j=1}^m \partial_{n+j+1}L(t)D_b^{\beta_j}\eta(t)+$$
 \begin{equation}\label{28}
 \sum_{p=0}^k  \partial_{m+n+p+2}L(t)\eta^{(p)}(t)]dt=0
 \end{equation}

By using the integration by parts formulas in the mentioned above
Lemma \ref{l1}, Lemma \ref{l2} and Lemma  \ref{l3} and the usual integration by parts formula , one obtains the following

$$\int_a^{b-\tau}[\sum_{i=1}^n D_{b-\tau}^{\alpha_i}(\partial_{i+1}L)(t)+\sum_{j=1}^m ~_aD^{\beta_j}(\partial_{n+j+1}L)(t)+$$
$$\sum_{p=0}^k(-1)^p\frac{d^p}{dt^p}~(\partial_{m+n+p+2}L)(t)+\sum_{p=0}^k(-1)^p\frac{d^p}{dt^p}~(\partial_{m+n+k+p+3}L)(t+\tau)-$$
$$\sum_{i=1}^n \frac{1}{\Gamma(\alpha_i)}D_{b-\tau}^{\alpha_i}(\int_{b-\tau}^b(D_b^{\alpha_i}\partial_{i+1}L)(s)(s-t)^{\alpha_i-1}ds)]\eta(t)dt+$$
$$\int_{b-\tau}^b[\sum_{i=1}^n D_b^{\alpha_i}(\partial_{i+1}L)(t)+\sum_{j=1}^m ~_{b-\tau}D^{\beta_j}(\partial_{n+j+1}L)(t)-$$
$$\sum_{j=1}^m\frac{1}{\Gamma(\beta_j)}~_{b-\tau}D^{\beta_j}(\int_a^{b-\tau}(~_aD^{\beta_j}\partial_{n+j+1}L)(s)(t-s)^{\beta_j-1}ds)]\eta(t)dt+$$
$$\sum_{p=0}^k(-1)^p\frac{d^p}{dt^p}~(\partial_{m+n+p+2}L)(t)]\eta(t)dt+$$
$$\sum_{p=1}^k\sum_{q=0}^{p-1}(-1)^q\frac{d^q}{dt^q}~(\partial_{m+n+p+2}L)(t)\eta^{p-q-1}(t)|_a^{b-\tau}+$$
$$\sum_{p=1}^k\sum_{q=0}^{p-1}(-1)^q\frac{d^q}{dt^q}~(\partial_{m+n+k+3}L)(t+\tau)\eta^{p-q-1}(t)|_a^{b-\tau}+$$
\begin{equation}\label{29}
\sum_{p=1}^k\sum_{q=0}^{p-1}(-1)^q\frac{d^q}{dt^q}~(\partial_{m+n+p+2}L)(t)\eta^{p-q-1}(t)|_{b-\tau}^b=0
\end{equation}

In equation (\ref{29}) if one  chooses $\eta$ such that $\eta(a)=0$ and $\eta\equiv 0$ on $[b-\tau,b]$, one gets
$$\sum_{i=1}^n D_{b-\tau}^{\alpha_i}(\partial_{i+1}L)(t)+\sum_{j=1}^m ~_aD^{\beta_j}(\partial_{n+j+1}L)(t)+$$
$$\sum_{p=0}^k(-1)^p\frac{d^p}{dt^p}~(\partial_{m+n+p+2}L)(t)+\sum_{p=0}^k(-1)^p\frac{d^p}{dt^p}~(\partial_{m+n+k+p+3}L)(t+\tau)-$$
\begin{equation}\label{30}
\sum_{i=1}^n \frac{1}{\Gamma(\alpha_i)}D_{b-\tau}^{\alpha_i}(\int_{b-\tau}^b(D_b^{\alpha_i}\partial_{i+1}L)(s)(s-t)^{\alpha_i-1}ds)=0.
\end{equation}
In equation (\ref{29}) if one  chooses $\eta$ such that $\eta^{(l)}(b)=0$ and $\eta\equiv 0$ on $[a,b-\tau]$, one gets
$$\sum_{i=1}^n D_b^{\alpha_i}(\partial_{i+1}L)(t)+\sum_{j=1}^m ~_{b-\tau}D^{\beta_j}(\partial_{n+j+1}L)(t)-$$
$$\sum_{j=1}^m\frac{1}{\Gamma(\beta_j)}~_{b-\tau}D^{\beta_j}(\int_a^{b-\tau}(~_aD^{\beta_j}\partial_{n+j+1}L)(s)(t-s)^{\beta_j-1}ds)+$$
\begin{equation}\label{31}
\sum_{p=0}^k(-1)^p\frac{d^p}{dt^p}~(\partial_{m+n+p+2}L)(t)=0.
\end{equation}
Now since both integrals in (\ref{29}) are now zero, one gets
$$\sum_{p=1}^k\sum_{q=0}^{p-1}(-1)^q\frac{d^q}{dt^q}~(\partial_{m+n+p+2}L)(t)\eta^{p-q-1}(t)|_a^{b-\tau}+$$
$$\sum_{p=1}^k\sum_{q=0}^{p-1}(-1)^q\frac{d^q}{dt^q}~(\partial_{m+n+k+3}L)(t)\eta^{p-q-1}(t+\tau)|_a^{b-\tau}+$$
\begin{equation}\label{32}
\sum_{p=1}^k\sum_{q=0}^{p-1}(-1)^q\frac{d^q}{dt^q}~(\partial_{m+n+p+2}L)(t)\eta^{p-q-1}(t)|_{b-\tau}^b=0.
\end{equation}
Thus one can state the following theorem

\begin{theorem}\label{t1}

Let $J(y)$ be a functional of the form (\ref{23})
defined on a set of continuous functions $y(t)$ which have
continuous Caputo fractional  left order derivatives of
orders $\alpha_i$ and right derivative of order $\beta_j$ in $[a,b]$ and
satisfy the conditions in (\ref{24}.~
Let $L:[a-\tau,b]\times \mathbb{R}^{m+n+2k+2}\rightarrow \mathbb{R}$ have continuous first and second partial
derivatives with respect to all of its arguments.  Then the  necessary
conditions that $J(y)$ possesses a minimum at  $y(x)$ are the
Euler-Lagrange equations
$$\sum_{i=1}^{n}D_{b-\tau}^{\alpha_i}(\frac{\partial L}{\partial ~a^CD^{\alpha_i}y(t)})(t)+\sum_{j=1}^{m}~_aD^{\beta_j}(\frac{\partial L}{\partial ~^CD_b^{\beta_j}y(t)})(t)+$$
$$\sum_{p=0}^{k}(-1)^p\frac{d^p}{dt^p}(\frac{\partial L}{\partial y^{(p)}(t)})(t)+\sum_{p=0}^{k}(-1)^p\frac{d^p}{dt^p}(\frac{\partial L}{\partial y^{(p)}(t-\tau)})(t+\tau)-$$
\begin{equation}\label{33}
\sum_{i=1}^n\frac{1}{\Gamma(\alpha_i)}D_{b-\tau}^{\alpha_i}(\int_{b-\tau}^b(D_b^{\alpha_i}(\frac{\partial L}{\partial ~_a^CD^{\alpha_i}y(t)})(s)(s-t)^{\alpha_i-1}ds)=0
\end{equation}
\noindent for $a\le t\le b-\tau$,
$$\sum_{i=1}^{n}D_{b}^{\alpha_i}(\frac{\partial L}{\partial ~a^CD^{\alpha_i}y(t)})(t)+\sum_{j=1}^{m}D_{b-\tau}^{\beta_j}(\frac{\partial L}{\partial ~^CD_b^{\beta_j}y(t)})(t)+$$
$$\sum_{p=0}^{k}(-1)^p\frac{d^p}{dt^p}(\frac{\partial L}{\partial y^{(p)}(t)})(t)-$$
\begin{equation}\label{34}
\sum_{j=1}^m\frac{1}{\Gamma(\beta_j)}~_{b-\tau}D^{\beta_j}(\int_{a}^{b-\tau}(~_aD^{\beta_j}(\frac{\partial L}{\partial ~^CD_b^{\beta_j}y(t)})(s)(t-s)^{\beta_j-1}ds)=0
\end{equation}
\noindent for $b-\tau\le t\le b$, and the transversality conditions
$$\sum_{p=1}^k\sum_{q=0}^{p-1}(-1)^q\frac{d^q}{dt^q}~(\frac{\partial L}{\partial y^{(p)}(t)})(t)\eta^{p-q-1}(t)|_a^{b-\tau}+$$
$$\sum_{p=1}^k\sum_{q=0}^{p-1}(-1)^q\frac{d^q}{dt^q}~(\frac{\partial L}{\partial y^{(p)}(t-\tau)}(t+\tau)\eta^{p-q-1}(t+\tau)|_a^{b-\tau}+$$
\begin{equation}\label{35}
\sum_{p=1}^k\sum_{q=0}^{p-1}(-1)^q\frac{d^q}{dt^q}~(\frac{\partial L}{\partial y^{(p)}(t)})(t)\eta^{p-q-1}(t)|_{b-\tau}^b=0.
\end{equation}

\noindent for any admissible function $\eta$ satisfying
$\eta(t)\equiv 0~t\in[a-\tau,a]$,\\$\eta^{(l)}(b)=0,~l=0,1,2,...,k-1$.
\end{theorem}

Theorem \ref{t1} can be generalized as follows

\begin{theorem}\label{t2}
  Consider the functional of the form
$$ J(y_1,y_2,...,y_d)=\int_{a}^{b}L[t,~_a^CD^{\alpha_1}y_1(t),~_a^CD^{\alpha_2}y_1(t),...,~_a^CD^{\alpha_n}y_1(t),$$
$$ ~_a^CD^{\alpha_1}y_2(t),~_a^CD^{\alpha_2}y_2(t),...,~_a^CD^{\alpha_n}y_2(t),...,$$
$$~_a^CD^{\alpha_1}y_d(t),~_a^CD^{\alpha_2}y_d(t),...,~_a^CD^{\alpha_n}y_d(t),$$
$$ ~^CD_b^{\beta_1}y_1(t),~^CD_b^{\beta_2}y_1(t),...,~^CD_b^{\beta_m}y_1(t),$$
$$ ~^CD_b^{\beta_1}y_2(t),~^CD_b^{\beta_2}y_2(t),...,~^CD_b^{\beta_m}y_2(t),...,$$
$$ ~^CD_b^{\beta_1}y_d(t),~^CD_b^{\beta_2}y_d(t),...,~^CD_b^{\beta_m}y_d(t),$$
$$ y_1(t),y_1'(t),...,y_1^{(k)}(t),y_2(t),y_2'(t),...,y_2^{(k)}(t),...,$$
$$ y_d(t),y_d'(t),...,y_d^{(k)}(t),y_1(t-\tau),y_1'(t-\tau),...,y_1^{(k)}(t-\tau),$$
\begin{equation}\label{36}
y_2(t-\tau),y_2'(t-\tau),...,y_2^{(k)}(t-\tau),...,y_d(t-\tau),y_d'(t-\tau),...,y_d^{(k)}(t-\tau)]dt,
 \end{equation}
\noindent  ,defined on sets of continuous functions
$y_i(x),~i=1,2,...,d$ that have left  Caputo fractional
derivatives of order $\alpha_i\in\mathbb{R},~i=1,2,...,n$  and right  Caputo
fractional derivatives of order $\beta_j\in\mathbb{R},~j=1,2,...m$ in the interval $[a,b]$
and satisfy the conditions
$$  y_i^{(l)}(b)=c_{il},~l=0,1,...,k ,~ y_i(t)=\phi_i(t) ~ i=1,2...,d,~ t\in[a-\tau,a]$$
\begin{equation}\label{37}
   a < b,~ \tau >0,~~~ \tau< b-a~,
\end{equation}
\noindent where $k-1\le\alpha_{\mbox{max}}<k, ~~\alpha_{\mbox{max}}=\mbox{max}\{\alpha_i,\beta_j\}_{1\le i\le n,1\le j\le m},$~$c_{il}$'s are constant and $\phi_i$'s are smooth
functions and  \\$L:[a-\tau,b]\times \mathbb{R}^{d(m+n+2k+2)}\rightarrow \mathbb{R}$ is a function with
continuous first and second partial derivatives with respect to all of its
arguments. For $y_i(x),~i=1,2,...,d$, satisfying
(\ref{37}) to be a minimum of (\ref{36}), it is
necessary  that
$$\sum_{i=1}^{n}D_{b-\tau}^{\alpha_i}(\frac{\partial L}{\partial ~a^CD^{\alpha_i}y_z(t)})(t)+\sum_{j=1}^{m}~_aD^{\beta_j}(\frac{\partial L}{\partial ~^CD_b^{\beta_j}y_z(t)})(t)+$$
$$\sum_{p=0}^{k}(-1)^p\frac{d^p}{dt^p}(\frac{\partial L}{\partial y_z^{(p)}(t)})(t)+\sum_{p=0}^{k}(-1)^p\frac{d^p}{dt^p}(\frac{\partial L}{\partial y_z^{(p)}(t-\tau)})(t+\tau)-$$
\begin{equation}\label{38}
\sum_{i=1}^n\frac{1}{\Gamma(\alpha_i)}D_{b-\tau}^{\alpha_i}(\int_{b-\tau}^b(D_b^{\alpha_i}(\frac{\partial L}{\partial ~_a^CD^{\alpha_i}y_z(t)})(s)(s-t)^{\alpha_i-1}ds)=0
\end{equation}
\noindent for $a\le t\le b-\tau$, $z=1,2,...,d$
$$\sum_{i=1}^{n}D_{b}^{\alpha_i}(\frac{\partial L}{\partial ~a^CD^{\alpha_i}y_z(t)})(t)+\sum_{j=1}^{m}D_{b-\tau}^{\beta_j}(\frac{\partial L}{\partial ~^CD_b^{\beta_j}y_z(t)})(t)+$$
$$\sum_{p=0}^{k}(-1)^p\frac{d^p}{dt^p}(\frac{\partial L}{\partial y_z^{(p)}(t)})(t)-$$
\begin{equation}\label{39}
\sum_{j=1}^m\frac{1}{\Gamma(\beta_j)}~_{b-\tau}D^{\beta_j}(\int_{a}^{b-\tau}(~_aD^{\beta_j}(\frac{\partial L}{\partial ~^CD_b^{\beta_j}y_z(t)})(s)(t-s)^{\beta_j-1}ds)=0
\end{equation}
\noindent for $b-\tau\le t\le b$,$~z=1,2,...,d$ and the transversality conditions
$$\sum_{p=1}^k\sum_{q=0}^{p-1}(-1)^q\frac{d^q}{dt^q}~(\frac{\partial L}{\partial y_z^{(p)}(t)})(t)\eta_z^{p-q-1}(t)|_a^{b-\tau}+$$
$$\sum_{p=1}^k\sum_{q=0}^{p-1}(-1)^q\frac{d^q}{dt^q}~(\frac{\partial L}{\partial y_z^{(p)}(t-\tau)}(t+\tau)\eta_z^{p-q-1}(t+\tau)|_a^{b-\tau}+$$
\begin{equation}\label{40}
\sum_{p=1}^k\sum_{q=0}^{p-1}(-1)^q\frac{d^q}{dt^q}~(\frac{\partial L}{\partial y_z^{(p)}(t)})(t)\eta_z^{p-q-1}(t)|_{b-\tau}^b=0.
\end{equation}

for any admissible vector function
$\eta=(\eta_1,\eta_2,...,\eta_d)$ satisfying\\
$\eta(t)\equiv (0,0,...,0)~t\in[a-\tau,a],~\eta^{(l)}(b)=(0,0,...,0),~d=0,1,2,...,k-1$.
\end{theorem}

\section{ The Fractional Optimal Control Problem}
 Find the optimal control variable $u(t)$ which minimizes the
 performance index
 $$J(y,u)=\int_{a}^{b}F[t,u(t),~_a^CD^{\alpha_1}y(t),~_a^CD^{\alpha_2}y(t),...,~_a^CD^{\alpha_n}y(t),$$
$$~^CD_b^{\beta_1}y(t),~^CD_b^{\beta_2}y(t),...
,~^CD_b^{\beta_m}y(t),y(t),y'(t),...,y^{(k)}(t),$$
\begin{eqnarray}\label{41}
y(t-\tau),y'(t-\tau),...,y^{(k)}(t-\tau)]dt,
\end{eqnarray}
\noindent subject to the constraint
$$G[t,u(t),~_a^CD^{\alpha_1}y(t),~_a^CD^{\alpha_2}y(t),...,~_a^CD^{\alpha_n}y(t),$$
$$~^CD_b^{\beta_1}y(t),~^CD_b^{\beta_2}y(t),...
,~^CD_b^{\beta_m}y(t),y(t),y'(t),...,y^{(k)}(t),$$
\begin{equation}\label{42}
y(t-\tau),y'(t-\tau),...,y^{(k)}(t-\tau)]=0
 \end{equation}
such that
$$y^{(l)}(b)=c_l,~l=0,1,2,...,k-1,~ y(t)=\phi(t)~ t\in [a-\tau,a],~a<b, \tau>0$$
\begin{equation}\label{43}
 \tau< b-a,~\alpha_i\in \mathbb{R},~i=1,2,...,n,~\beta_j\in\mathbb{R},~j=1,2,...,m,
\end{equation}
\noindent where $c_{l}$ are constant and $F$ and $G$ are  functions  $ [a-\tau,b]\times\mathbb{R}^{n+m+2k+3}\rightarrow  \mathbb{R}$
with continuous first and second partial derivatives with respect
to all of their arguments arguments.

 To find the optimal control, one defines a modified performance
 index as
\begin{equation}\label{44}
 \hat{J}(y,u)=\int_{a}^{b}
F+\lambda(t)G~dt,
 \end{equation}
\noindent where $\lambda$ is a Lagrange multiplier  or an adjoint
variable. Using the conditions (\ref{38}), (\ref{39})
and (\ref{40}) in Theorem \ref{t2}, the following necessary equations for optimal control are found:
Euler-Lagrange equations
$$\sum_{i=1}^{n}D_{b-\tau}^{\alpha_i}(\frac{\partial F}{\partial ~a^CD^{\alpha_i}y(t)})(t)+
\sum_{i=1}^{n}D_{b-\tau}^{\alpha_i}(\lambda\frac{\partial G}{\partial ~a^CD^{\alpha_i}y(t)})(t)+$$
$$\sum_{j=1}^{m}~_aD^{\beta_j}(\frac{\partial F}{\partial ~^CD_b^{\beta_j}y(t)})(t)+       \sum_{j=1}^{m}~_aD^{\beta_j}(\lambda\frac{\partial G}{\partial ~^CD_b^{\beta_j}y(t)})(t)+$$
$$\sum_{p=0}^{k}(-1)^p\frac{d^p}{dt^p}(\frac{\partial F}{\partial y^{(p)}(t)})(t)+
\sum_{p=0}^{k}(-1)^p\frac{d^p}{dt^p}(\lambda\frac{\partial G}{\partial y^{(p)}(t)})(t)+$$

$$+\sum_{p=0}^{k}(-1)^p\frac{d^p}{dt^p}(\frac{\partial F}{\partial y^{(p)}(t-\tau)})(t+\tau)++\sum_{p=0}^{k}(-1)^p\frac{d^p}{dt^p}(\lambda\frac{\partial G}{\partial y^{(p)}(t-\tau)})(t+\tau)-$$
$$\sum_{i=1}^n\frac{1}{\Gamma(\alpha_i)}D_{b-\tau}^{\alpha_i}(\int_{b-\tau}^b(D_b^{\alpha_i}(\frac{\partial F}{\partial ~_a^CD^{\alpha_i}y(t)})(s)(s-t)^{\alpha_i-1}ds)-$$
$$\sum_{i=1}^n\frac{1}{\Gamma(\alpha_i)}D_{b-\tau}^{\alpha_i}(\int_{b-\tau}^b(D_b^{\alpha_i}(\lambda\frac{\partial G}{\partial ~_a^CD^{\alpha_i}y(t)})(s)(s-t)^{\alpha_i-1}ds)+$$
\begin{equation}\label{45}
\frac{\partial F}{\partial u(t)}(t)+\lambda(t)\frac{\partial G}{\partial u(t)}(t)=0,
\end{equation}
\noindent for $a\le t\le b-\tau$,
$$\sum_{i=1}^{n}D_{b}^{\alpha_i}(\frac{\partial F}{\partial ~a^CD^{\alpha_i}y(t)})(t)+
\sum_{i=1}^{n}D_{b}^{\alpha_i}(\lambda\frac{\partial G}{\partial ~a^CD^{\alpha_i}y(t)})(t)+$$

$$\sum_{j=1}^{m}D_{b-\tau}^{\beta_j}(\frac{\partial F}{\partial ~^CD_b^{\beta_j}y(t)})(t)+
\sum_{j=1}^{m}D_{b-\tau}^{\beta_j}(\lambda\frac{\partial G}{\partial ~^CD_b^{\beta_j}y(t)})(t)+$$
$$\sum_{p=0}^{k}(-1)^p\frac{d^p}{dt^p}(\frac{\partial F}{\partial y^{(p)}(t)})(t)+\sum_{p=0}^{k}(-1)^p\frac{d^p}{dt^p}(\lambda\frac{\partial G}{\partial y^{(p)}(t)})(t)+$$
$$\sum_{j=1}^m\frac{1}{\Gamma(\beta_j)}~_{b-\tau}D^{\beta_j}(\int_{a}^{b-\tau}(~_aD^{\beta_j}(\frac{\partial F}{\partial ~^CD_b^{\beta_j}y(t)})(s)(t-s)^{\beta_j-1}ds)+$$
$$\sum_{j=1}^m\frac{1}{\Gamma(\beta_j)}~_{b-\tau}D^{\beta_j}(\int_{a}^{b-\tau}(~_aD^{\beta_j}(\lambda\frac{\partial L}{\partial ~^CD_b^{\beta_j}y(t)})(s)(t-s)^{\beta_j-1}ds)+$$
\begin{equation}\label{46}
\frac{\partial F}{\partial u(t)}(t)+\lambda(t)\frac{\partial G}{\partial u(t)}(t)=0,
\end{equation}
\noindent for $b-\tau\le t\le b$, and the transversality conditions
$$\sum_{p=1}^k\sum_{q=0}^{p-1}(-1)^q\frac{d^q}{dt^q}~(\frac{\partial F}{\partial y^{(p)}(t)})(t)\eta^{p-q-1}(t)|_a^{b-\tau}+$$
$$\sum_{p=1}^k\sum_{q=0}^{p-1}(-1)^q\frac{d^q}{dt^q}~(\lambda\frac{\partial G}{\partial y^{(p)}(t)})(t)\eta^{p-q-1}(t)|_a^{b-\tau}+$$
$$\sum_{p=1}^k\sum_{q=0}^{p-1}(-1)^q\frac{d^q}{dt^q}~(\frac{\partial F}{\partial y^{(p)}(t-\tau)}(t+\tau)\eta^{p-q-1}(t+\tau)|_a^{b-\tau}+$$
$$\sum_{p=1}^k\sum_{q=0}^{p-1}(-1)^q\frac{d^q}{dt^q}~(\lambda\frac{\partial G}{\partial y^{(p)}(t-\tau)}(t+\tau)\eta^{p-q-1}(t+\tau)|_a^{b-\tau}+$$
$$\sum_{p=1}^k\sum_{q=0}^{p-1}(-1)^q\frac{d^q}{dt^q}~(\frac{\partial F}{\partial y^{(p)}(t)})(t)\eta^{p-q-1}(t)|_{b-\tau}^b$$
\begin{equation}\label{47}
\sum_{p=1}^k\sum_{q=0}^{p-1}(-1)^q\frac{d^q}{dt^q}~(\lambda\frac{\partial G}{\partial y^{(p)}(t)})(t)\eta^{p-q-1}(t)|_{b-\tau}^b=0,
\end{equation}
\noindent where $\eta$ is any admissible function satisfying
$\eta(t)\equiv 0~t\in[a-\tau,a]$,\\$\eta^{(l)}(b)=0,~l=0,1,2,...,k-1$.

\section{Conclusion}
In this manuscript we have developed a fractional control problem
in the presence of both left and right Caputo fractional derivatives of any order
 and delay in the state variables and their derivatives. The
results were applied in order to find the necessary conditions for
the optimal control problems. When $\alpha_i\rightarrow 1,~\forall i=1,...,n$ and $~\beta_j \rightarrow 1,~\forall j=1,...,m$ the classical problem is recovered.

\section*{Acknowledgments}
This work is partially supported by the Scientific and Technical
Research Council of Turkey.


\begin{thebibliography}{99}
\bibitem{trujillo}
A. A. Kilbas  , H. H. Srivastava  and  J. J. Trujillo , { \em Theory
and Applications of Fractional Differential Equations\/},
Elsevier, Amsterdam 2006

\bibitem{samko} S. G. Samko , A. A. Kilbas and O.I. Marichev,
 {\em Fractional Integrals and Derivatives - Theory and Applications\/}
 Gordon and Breach, Linghorne, P.A. 1993.

 \bibitem{podlubny1}
 I. Podlubny,  {\em Fractional Differential Equations\/}, Academic
Press, San Diego CA 1999.

\bibitem{richard}
R. L. Magin, {\em  Fractional Calculus in Bioengineering\/},
Begell House Publisher, Inc. Connecticut 2006.

 \bibitem{west}  B. J. West, M.  Bologna, P. Grigolini,
~\emph{\textit{}Physics of Fractal operators},~Springer,~New
York 2003.


\bibitem{podlubny}
N. Heymans, I. Podlubny, Physical interpretation of
initial conditions for fractional differential equations with
Riemann-Liouville fractional derivatives,~{\em Rheologica Acta}
{\bf 45} (2006) 765-771.


\bibitem{machado}
I. S. Jesus, J. A. T.  Machado,  Fractional control of
heat diffusion system, \emph{Nonlinear Dyn.}, \textbf{54(3)
} (2008) 263-282.


\bibitem{mainardi}
F. Mainardi, Y. Luchko, G.  Pagnini, The
fundamental solution of the space-time fractional diffusion
equation {\em Frac. Calc. Appl. Analys.\/} {\bf
4(2)} (2001)  153-192.


\bibitem{thab}
B. Dumitru , T. Maraaba, F. Jarad, Fractional Principles with
Delay, {\em J. Phys. A: Math. Theor.} {\bf 41}(2008)
doi:10.1088/1751-8113/41/31/315403.

\bibitem{fahd}
F. Jarad, T. Maraaba , D. Baleanu  Fractional Variational Principles with
Delay  within Caputo Derivatives, {\em Rep. Math. Phys.} {\bf 65(1)} (2010)  17-28.


\bibitem{momani}
S. Momani,  A numerical scheme for the solution of multi-order
fractional differential equations,  {\em   Appl. Math. Compt.} {\bf
182} (2006)  761-786.

\bibitem{mainardi}
F. Mainardi , Y. Luchko, G. Pagnini,  The fundamental
solution of the space-time fractional the space-time fractional
diffusion equation {\em Frac. Calc. Appl. Analys.\/} {\bf 4(2)}(2001)
153-192.

\bibitem{enrico}
 E. Scalas ,  R. Gorenflo,  F. Mainardi, Uncoupled continuous-time
random walks: Solution and limiting behavior of the master
equation, {\em Phys. Rev. E\/} \textbf{69} (2004) 011107-1.


\bibitem{podlubny1}
N. Heymans,  I.  Podlubny,  Physical interpretation of initial
conditions for fractional differential equations with
Riemann-Liouville fractional derivatives, {\em Rheol. Acta} {\bf
45}2006 765-771.

\bibitem{riewe1} F. Riewe, Nonconservative Lagrangian and Hamiltonian
mechanics,  {\em  Phys. Rev. E\/} {\bf 53}(1996) 1890-1899.
\bibitem{riewe2} F. Riewe, ~Mechanics with fractional derivatives, {\em
Phys. Rev. E\/} {\bf 55}(1997) 3581-3592.
\bibitem{klimek1} M. Klimek,  Fractional sequential mechanics - models
with symmetric fractional derivative, {\em  Czech. J. Phys.\/} {\bf
51}(2001) 1348-1354.
\bibitem{klimek2}
M. Klimek,  Lagrangean and Hamiltonian fractional sequential
mechanics, {\em Czech. J. Phys.\/} {\bf 52}(2002) 1247-1253

\bibitem{agrawal5}
D. Baleanu, O. P.  Agrawal,  Fractional Hamilton formalism
within Caputo's derivative,~
 {\em Czech. J. Phys.\/} {\bf 56}(2006)
1087-1092.

\bibitem{agrawal1}O. P. Agrawal,
  Formulation of Euler-Lagrange equations for fractional variational
problems,  {\em J. Math. Anal. Appl.\/} {\bf 272} (2002)2006  368-379.
\bibitem{agrawal2}O. P. Agrawal,  Fractional variational calculus and the
transversality conditions,  {\em J. Phys.A.: Math.Gen.\/} {\bf 39}(2006)
10375-10384
 \bibitem{agrawal3}
 O. P. Agrawal, Generalized Euler-Lagrange equations and transversality
conditions for FVPs in terms of the Caputo
 derivative,
 \emph{J. Vib. Contr.} {\bf 13(9-10)}(2007) 1217-1237.

\bibitem{agrawal4}
O. P. Agrawal, D. Baleanu, Hamiltonian formulation and a
direct numerical scheme for Fractional Optimal Control Problems,
 {\em J. Vib. Contr.\/} {\bf 13(9-10)}(2007) 1269-1281.

 \bibitem{chen}Y. Q. Chen ,B. M. ~Vinagre, I. Podlubny,
Continued fraction expansion approaches to discretizing
fractional order derivatives-an expository review, \emph{Nonlinear
Dyn. }\textbf{38(1-4)} (2004) 155-170.
\bibitem{rosen} J. F. Rosenblueth, Systems with time delay in the
calculus of variations : A variational approach, \emph{ J. Math.
Contr. Inf.}\textbf{(5)} (1988) 125-145.

\bibitem{zaslavsky}
V. E. Tarasov, G. M. Zaslavsky,  Nonholonomic constraints with
fractional derivatives, \emph{J. Phys. A-Math. Gen.
}\textbf{39(31)}(2006) 9797-9815.
\bibitem{agragen} O. P. Agrawal,  A general formulation and
solution scheme for fractional optimal control problems, \emph
{Nonlinear Dyn.}\textbf{38}(2004) 323-337.

  \bibitem{baleanu1}
 E. M. Rabei, K. I. Nawafleh, R. S.  Hijjawi, S. I. Muslih,~
D. Baleanu, The Hamilton formalism with fractional derivatives,
{\em J. Math. Anal. Appl.\/} {\bf 327}(2007)  891-897.

\bibitem{baleanu2}
 D. Baleanu, S. I.   Muslih ,  Lagrangian formulation of classical fields
within Riemann-Liouville fractional derivatives, {\em Physica
Scripta\/} {\bf 72(2-3)} (2005) 119-121.

\bibitem{baleanu3}
  D. Baleanu, S. I.   Muslih, Hamiltonian formulation of systems with
linear velocities within Riemann-Liouville fractional derivatives,
 {\em J. Math. Anal. Appl.\/} {\bf
 304(3)}(2005)  599-606.

\bibitem{baleanu4}
 D. Baleanu, T. Avkar,   Lagrangians with linear velocities within
Riemann-Liouville fractional derivatives,
 {\em Nuovo Cimento B\/}
{\bf 119}(2004)  73-79

\bibitem{driver}  R. D. Driver,  Ordinary and Delay Differential Equations, Springer-Verlag, New York, 1977.

\bibitem{chan}
W. Deng, C. Li,J. Lu,  Stability analysis of linear
fractional differential system with multiple time scales,
\emph{Nonlinear Dyn.} \textbf{48}(2007)  409-416.

\bibitem{bliss}
G. A. Bliss,  Lectures on the Calculus of variations,
The University of Chicago Press, Chicago and London, 1963.

\bibitem{gregory}
J. Gregory, C. Lin, Constrained optimization n the
calculus of variations and optimal control theory, Van Nostrand
Reinhold, 1989.

\bibitem{Littlewood} G. H. Hardy, J. E.  Littlewood, Some properties
of fractional integrals, \emph{Math. Zeitschrift}~\textbf{27} (1928) 565-606.
\end{thebibliography}
\end{document}